\documentstyle[twoside,amssymb]{article}

\newtheorem{theorem}{Theorem}[section]
\newtheorem{corollary}[theorem]{Corollary}
\newtheorem{conjecture}[theorem]{Conjecture}
\newtheorem{lemma}[theorem]{Lemma}

\input tcilatex
\QQQ{Language}{
American English
}
\begin{document}

\author{Alexander Dvorsky and\ Siddhartha Sahi\thanks{
\ \thinspace This work was supported by an NSF grant.}\medskip \\
Department of Mathematics, Rutgers University, New Brunswick\\
NJ 08903, USA}
\title{Tensor Products of Singular Representations and an Extension of the $\theta $%
-correspondence\bigskip\ }
\date{}
\maketitle

\begin{abstract}
In this paper we consider the problem of decomposing tensor products of
certain singular unitary representations of a semisimple Lie group $G$.
Using explicit models for these representations (constructed earlier by one
of us) we show that the decomposition is controlled by a reductive
homogeneous space $G^{\prime }/H^{\prime }$. Our procedure establishes a
correspondence between certain unitary representations of $G$ and those of $%
G^{\prime }$. This extends the usual $\theta $--correspondence for dual
reductive pairs. As a special case we obtain a correspondence between
certain representations of real forms of $E_7$ and $F_4$.
\end{abstract}

\date{\ \ }

\section{Introduction}

Let $F$ be a field and $\varepsilon $ some fixed additive character of $F$.
If $W$ is a finite dimensional vector space over $F$ endowed with a
non-degenerate skew-symmetric form, we can consider an associated Heisenberg
group $H(W)$. Denote by $\rho _\varepsilon $ an irreducible unitary
representation of $H(W)$ on which the center of $H(W)$ operates via the
character $\varepsilon $ (it is unique by the theorem of Stone and von
Neumann). Since the symplectic group $Sp(W)$ operates on $H(W)$ via its
action on the vector space $W$, it also acts on the representation $\rho
_\varepsilon $. The action is trivial on the center of $H(W)$ and therefore,
for any $g\in Sp(W),$ there is an operator $\omega _\varepsilon (g)$ (unique
up to scaling) which intertwines $\rho _\varepsilon $ with $\rho
_\varepsilon ^g$. These operators form the {\em oscillator representation }%
of $Sp(W)$. In general, the oscillator representation $\omega _\varepsilon $
is projective, but it always corresponds to an ordinary representation of a
two-fold cover of $Sp(W)$. This cover is denoted by $\widetilde{Sp(W)}$ and
is called the {\em metaplectic group}.

The oscillator representation originated in the works of Segal and Shale and
was immortalized by Weil, who used it in his construction of the
theta-functions on the metaplectic group \cite{weil}. In one of the earliest
works on the spectrum of the Weil (oscillator) representation, Gelbart \cite
{gelbart} studied the decomposition of the tensor product $\omega ^{\prime
}=\omega ^{\otimes k}$, where $\omega $ is the oscillator representation of
the real symplectic group $Sp(2m,{\Bbb R})$ and $k\geq 2m$. He demonstrated
that for $k=2m$ all representations of the holomorphic discrete series for $%
Sp(2m,{\Bbb R})$ occur in the spectrum of $\omega ^{\prime }$. Kashiwara and
Vergne \cite{kashiwara-vergne} extended the results of \cite{gelbart} to
tensor products $\omega ^{\otimes k}$, $k\geq 1$. In particular, any unitary
highest weight representation of $Sp(2m,{\Bbb R})$ will appear in the
decomposition of $\omega ^{\otimes k}$ for some appropriate value of $k$.

Later, the approach of \cite{gelbart} and \cite{kashiwara-vergne} was
replaced by the modern approach to the $\theta $-correspondence. One starts
with a reductive dual pair of algebraic groups $G$ and $G^{\prime }$,
defined over some local field $F$, which are mutual centralizers inside a
symplectic group $Sp$. Let $\omega $ be an oscillator representation of the
metaplectic group $\widetilde{Sp}$ on a Hilbert space ${\cal H}$, and let $%
\overline{E}\subset $ $\widetilde{Sp}$ denote the preimage of a reductive
subgroup $E\subset Sp$. Denote by ${\cal R}(E)$ the set of (equivalence
classes) of continuous irreducible representations of $\overline{E}$ on a
locally convex space, which can be realized as quotients of ${\cal H}^\infty 
$ by $\omega (\overline{E})$--invariant subspaces.

\begin{conjecture}[Howe's duality conjecture]
The set ${\cal R}(G\cdot G^{\prime })$ is the graph of a bijection between
all of ${\cal R}(G)$ and all of ${\cal R}(G^{\prime })$. Moreover, an
element of ${\cal R}(G\cdot G^{\prime })$ occurs as a quotient of $\omega $
in a unique way.
\end{conjecture}

Howe's conjecture has been proved for $F={\Bbb R}$ or ${\Bbb C}$ \cite
{howe-tr}, and also for all non-Archimedean local fields of odd residue
characteristic. The resulting correspondence between the irreducible
representations of $\overline{G}$ and those of $\overline{G^{\prime }}$ is
called the Howe duality correspondence, or the $\theta $-correspondence. In
general, this correspondence does not preserve unitarity, i.e., a unitary
representation of $\overline{G}$ can correspond to a non-unitarizable
representation of $\overline{G^{\prime }}$.

If one member of a dual pair (say $G^{\prime }$) is much ``smaller'' than
the second group, unitarity is preserved and the duality correspondence is a
particularly nice one (this is a so called stable range duality \cite
{howe-rank}, \cite{li-singular}). We denote by $\widehat{G}(\varepsilon )$
and $\widehat{G^{\prime }}(\varepsilon )$ the subsets of the unitary duals
of $\overline{G}$ and $\overline{G}^{\prime }$ consisting of those unitary
irreducible representations whose restriction to the kernel of the
projection $\widetilde{Sp}\rightarrow Sp$ (i.e., the group ${\Bbb Z}_2$) is
a multiple of the non-trivial character of ${\Bbb Z}_2$. Then the Howe
correspondence gives an injection $\widehat{G^{\prime }}(\varepsilon )$ $%
\hookrightarrow $ $\widehat{G}(\varepsilon )$. In many cases the coverings $%
\overline{G}\rightarrow G$ and $\overline{G}^{\prime }\rightarrow G^{\prime
} $ are trivial, and we obtain an injection of the unitary dual of $%
G^{\prime } $ into that of $G$.\smallskip\ 

{\em Example}. $G^{\prime }=U(1)$ and $G=U(p,q)$ form a stable range dual
pair inside $Sp(2p+2q,{\Bbb R})$. The representations of $U(p,q)$ appearing
in the $\theta $-correspondence for this reductive dual pair are the
``ladder'' representations, which were given this name because their $%
\overline{K}$-types lie along a line, i.e., their highest weights are
obtained from the highest weight of the lowest $\overline{K}$-type by adding
multiples of a single vector.\smallskip\ 

The duality conjecture also has a global version \cite[2.5]{gelbart(crm)},
when $F$ is a global field, ${\Bbb A}$ is an adele ring of $F$, and one
considers a dual pair inside a symplectic group over ${\Bbb A}$. This
version is of considerable representation- and number-theoretic interest,
since the global duality correspondence provides a way to lift automorphic
forms between members of a dual pair \cite{gelbart(crm)}, \cite{prasad}. For
example, if $\overline{G^{\prime }}=\widetilde{SL(2)}$ and $G\simeq PGL(2)$
is realized as an orthogonal group preserving a 3-variables quadratic form $%
Q(x_1,x_2,x_3)=x_1^2-x_2x_3$, the $\theta $-correspondence produces the
Shimura lifting which associates a modular form of weight $n-1$ to a modular
form of half-integral weight $\frac n2$.

It is a remarkable fact that {\em every} classical group and {\em no}
exceptional group can be realized as a member of a reductive dual pair. It
seems desirable to determine whether one can extend this theory to other
reductive groups in a reasonable manner. In this paper we attempt to extend
the original (\cite{gelbart}, \cite{kashiwara-vergne}) approach to the $%
\theta $-correspondence. We study the decompositions of tensor products of
certain small unitary representations of a real reductive group $G$ and
construct a parametrization of the spectra of these tensor products.

Let $\Omega $ be a symmetric tube domain of rank $n$, and $G=\limfunc{Aut}%
(\Omega )$. The Shilov boundary of $\Omega $ is of the form $G/P$ where $%
P=LN $ is the Siegel--Poincar\'{e} parabolic subgroup of $G$. The nilradical 
$N$ is abelian and so is isomorphic to its Lie algebra ${\frak n}$, and the
Levi subgroup $L$ has finitely many (coadjoint) orbits on ${\frak n}^{*}$.
These orbits are indexed by their ``signatures'', where a signature $p$
consists of a pair of non-negative integers\thinspace $p=(p^{+},p^{-})$ with 
$\left| p\right| \stackrel{\text{def}}{=}p^{+}+p^{-}\leq n$.

Every non-open orbit (with $\left| p\right| <n)$ has a canonical $L$%
-equivariant measure, and the main result of \cite{sahi-expl} is the
construction of an irreducible representation of the universal cover of $G$
on the associated $L^2$-space. Usually these representations descend to $G$,
and in all cases they can be viewed as linear representations of a certain
double cover of $G$, denoted by $\overline{G}$. In this paper we consider
the problem of decomposing certain tensor products of these representations.
More precisely, let $\Pi =\pi _1\otimes \cdots \otimes \pi _s$ be the tensor
product of representations associated to orbits ${\cal O}_{p_1}$,$\ldots $,$%
{\cal O}_{p_s}$ whose signatures $p_i=(p_i^{+},p_i^{-})$ satisfy 
\begin{equation}
(p_1^{+}+p_1^{-})+\cdots +(p_s^{+}+p_s^{-})\leq n\text{.}
\label{-stable-range}
\end{equation}

Under assumption (\ref{-stable-range}), which is an analogue of the stable
range condition for the $\theta $-correspondence, each of the following
spaces contains an open and dense $L$-orbit 
\[
{\cal O}\equiv {\cal O}_{p_1}+\ldots +{\cal O}_{p_s},\;\;{\cal O}^{\prime
}\equiv {\cal O}_{p_1}\times \ldots \times {\cal O}_{p_s}\text{.} 
\]
Fix a generic point $\xi ^{\prime }=(\xi _1,\ldots ,\xi _k)$ in ${\cal O}%
^{\prime }$ such that $\xi =\xi _1+\ldots +\xi _k$ is generic in ${\cal O}$.
We denote the inverse images of $L$ and $P$ in $\overline{G}$ by $\overline{L%
}$ and $\overline{P}$, respectively. Let $S$ and $S^{\prime }$ be the
stabilizers of $\xi ^{\prime }$ and $\xi $ in $\overline{L}$, and let $\chi
_\xi $ be the unitary character of $N$ defined by $\chi _\xi (\exp x)%
\stackrel{\text{def}}{=}\exp (i\left\langle \xi ,x\right\rangle )$ for $x\in 
{\frak n}$.

\begin{theorem}
The restriction of $\Pi $ to $\overline{P}$ is isomorphic to $\nu \otimes 
\limfunc{Ind}_{S^{\prime }N}^{\overline{P}}(1\otimes \chi _\xi )$, where $%
\nu $ is a certain unitary character of $\overline{P}$.\label{-th00}
\end{theorem}

In general, $S,S^{\prime }$ are not reductive, however they contain
reductive groups $G^{\prime }$ and $H^{\prime }$ and a common normal
subgroup $N$ such that $S=G^{\prime }\ltimes Z$ and $S^{\prime }=H^{\prime
}\ltimes Z$. We consider the direct integral decomposition 
\[
L^2(G^{\prime }/H^{\prime })=\int_{\pi \in \widehat{G^{\prime }}}^{\oplus
}m(\pi )\pi \,d\mu (\pi )\text{.} 
\]

For each $\pi $ occurring in $L^2(G^{\prime }/H^{\prime }),$ we define $%
\Theta (\pi )\in \widehat{\overline{P}}$ by 
\[
\Theta (\pi )=\nu \otimes \limfunc{Ind}\nolimits_{G^{\prime }ZN}^{\overline{P%
}}(\pi \otimes 1\otimes \chi _\xi )\text{.} 
\]
By Mackey theory, $\Theta (\pi )$ is irreducible, and Theorem \ref{-th00}
implies

\begin{theorem}
The restriction of $\Pi $ to $\overline{P}$ has the decomposition 
\begin{equation}
\Pi |_{\overline{P}}=\int_{\pi \in \widehat{G^{\prime }}}^{\oplus }m(\pi
)\Theta (\pi )\,d\mu (\pi )\text{ .}  \label{-pi-dirint}
\end{equation}
\label{-th01}
\end{theorem}

The main result of the paper is the following

\begin{theorem}
For almost every $\pi $ (with respect to the Plancherel measure $d\mu $), $%
\Theta (\pi )$ extends to an irreducible representation of $\overline{G}$,
so that {\em (\ref{-pi-dirint})} is also a $\overline{G}$-decomposition. If $%
\sum_{i=1}^s|p_i|<n$, this extension is unique.\label{-g-decomp}
\end{theorem}

Thus the map $\pi \rightarrow \Theta (\pi )$ gives a (measurable) bijection
between unitary representations of $\overline{G}$ occurring in $\Pi $ and
the unitary representations of $G^{\prime }$ occurring in $L^2(G^{\prime
}/H^{\prime })$. \ 

We now discuss some special cases of the above result:

If $s=2$, then $G^{\prime }/H^{\prime }$ is a symmetric space, which is
Riemannian if and only if ${\cal O}_{p_1}$ and ${\cal O}_{p_2}$ both have
definite signatures (of the form $(p^{+},0)$ or $(0,p^{-})$). Positive
(resp. negative) definite orbits correspond to the highest (resp. lowest)
weight singular representations of $\overline{G}$, and in this case our
constructions complement the results on the tensor products of holomorphic
and anti-holomorphic discrete series representations in \cite{repka}.

For Riemannian symmetric spaces, and also for several non-Riemannian ones,
we have $m(\pi )\leq 1$. Thus in these cases we deduce that $\pi _1\otimes
\pi _2$ is multiplicity free.

If $\Omega $ is the Siegel upper half plane then $G=Sp(2n,{\Bbb R})$, $%
\overline{G}$ is the metaplectic group and $G^{\prime }/H^{\prime
}=O(p^{+},p^{-})/[O(p_1^{+},p_1^{-})\times \ldots \times O(p_s^{+},p_s^{-})]$%
, where $p^{+}=\sum p_i^{+}$, $p^{-}=\sum p_i^{-}$. In this case our
correspondence coincides with the $H^{\prime }$-spherical part of the $%
\theta $-correspondence.

Finally, if $\Omega $ is the exceptional tube domain, then $\overline{G}=G$
is the simply connected exceptional group $E_{7(-25)}$. If we take $s=2$ and 
$p_1=(1,0)$, then among the possibilities for $G^{\prime }/H^{\prime }$ are
the various forms of the Cayley projective plane \cite[p.118]{adams}, i.e.,
for $p_2=(2,0),$ $p_2=(0,2)$ and $p_2=(1,1)$, we obtain respectively 
\[
F_{4(-52)}/Spin(9),\;F_{4(-20)}/Spin(9)\text{ and}\;F_{4(-20)}/Spin(1,8). 
\]
Note that these symmetric spaces are multiplicity-free (see \cite{vandijk}
for the non-Riemannian space $F_{4(-20)}/Spin(1,8)$) and by Theorem \ref
{-g-decomp}, so is $\Pi =\pi _1\otimes \pi _2$.

If we take $p_1=p_2=(1,0)$, then $\Pi $ is a tensor square of a highest
weight representation $\pi _1$, and a description of the spectrum of this
tensor square is a key step in the classification of unitarizable highest
weight modules in \cite{ehw}.

Just as with the $\theta $-correspondence, we expect that our results will
have smooth and global analogues. We shall take up some of these questions
in subsequent papers.

\section{Notation and Preliminaries}

\subsection{Groups and subgroups}

Let $G$ be one of the following groups:

\begin{itemize}
\item  $Sp(2n,{\Bbb R})$ (case I1),
\end{itemize}

\begin{itemize}
\item  $U(n,n)$ (case I2),
\end{itemize}

\begin{itemize}
\item  $O^{*}(4n)$ (case I3),
\end{itemize}

\begin{itemize}
\item  $O(2,j)$ (case \thinspace I4),
\end{itemize}

\begin{itemize}
\item  $E_{7(-25)}$ (case I5),
\end{itemize}

and $K$ be the maximal compact subgroup of $G$. Then $\Omega =G/K$ is a
symmetric domain of tube type \cite[p.474]{helgason}. Taking $G=U(n,n)$ or $%
O(2,j)$ instead of $SU(n,n)$ or $SO(2,j)$ is not really necessary, but will
make some arguments more straightforward.

The restricted root system for each of the groups listed above is of type $%
C_n$, where $n$ is the real rank of the group $G$. Let $\Delta =\{\beta
_1,\beta _2,\ldots ,\beta _n\}$ be the basis of the restricted root system,
enumerated in such a way that the corresponding Dynkin diagram is 
\[
\begin{array}{ccc}
\stackrel{\beta _1}{\circ }\text{------}\stackrel{\beta _2}{\circ }\text{%
------}\stackrel{\beta _3}{\circ } & \cdots & \stackrel{\beta _{n-2}}{\circ }%
\text{------}\stackrel{\beta _{n-1}}{\circ }\Longleftarrow \stackrel{\beta _n%
}{\circ }
\end{array}
\text{ .} 
\]
There exists a one-to-one correspondence between the set of maximal
parabolic subgroups of $G$ and the set of maximal subsets of $\Delta $. We
will be interested in two parabolic subgroups of $G$ -- the Siegel parabolic
of $G$ (this corresponds to the set $\Delta \setminus \{\beta _n\}$) and the
maximal parabolic subgroup corresponding to the set $\Delta \setminus
\{\beta _1\}$. We denote the first by $P$ and the second by \thinspace $%
P^{\prime }$.

The Levi decomposition of $P$ is $P=L\cdot N$, where the subgroup $N$ is
abelian (e.g., for $G=O(2,j)$ we have $n=2,$\thinspace $L={\Bbb R}^{*}\times
O(1,j-1)$ and $N={\Bbb R}^{1,j-1}$).

The Langlands decomposition of $P^{\prime }$ is $P^{\prime }=M^{\prime
}AN^{\prime }$, where the radical $N^{\prime }$ is a two-step nilpotent
group with a one-dimensional center{\bf \ }$ZN^{\prime }$, and we can
identify it with the real Heisenberg group of dimension $2m+1$. The vector
subgroup $A$ is one-dimensional, i.e., $A={\Bbb R}^{*}$. For example, $%
G=E_{7(-25)}$ gives $M^{\prime }A=SO(2,10)\times {\Bbb R}^{*}$ and $%
N^{\prime }$ is the Heisenberg group associated with a 32-dimensional real
vector space.

The group $M^{\prime }$ splits into a direct product of a compact factor and
a noncompact group, which we denote by $G_{-}$. The group $G_{-}$ belongs to
one of the families (I1)-(I4), and we can consider its Siegel parabolic $%
P_{-}$, the nilradical $N_{-}$ of $P_{-}$, etc. In general, all subgroups of 
$G_{-}$ will be written with a minus as a subscript.

The information about some of the subgroups we defined above is summarized
in the following table.\bigskip\ \newline\ {\small 
\begin{tabular}{|c|c|c|c|c|}
\hline
$G$ & $N$ & $M^{\prime }$ & $G_{-}$ & $m$ \\ \hline
\multicolumn{1}{|l|}{$Sp(2n,{\Bbb R)}$} & \multicolumn{1}{|l|}{$Sym(n,{\Bbb R%
})$} & \multicolumn{1}{|l|}{$Sp(2n-2,{\Bbb R)}$} & \multicolumn{1}{|l|}{$%
Sp(2n-2,{\Bbb R)}$} & $n-1$ \\ \hline
\multicolumn{1}{|l|}{$U(n,n)$} & \multicolumn{1}{|l|}{$Herm(n,{\Bbb C})$} & 
\multicolumn{1}{|l|}{$U(1)\times U(n-1,n-1)$} & \multicolumn{1}{|l|}{$%
U(n-1,n-1)$} & $2(n-1)$ \\ \hline
\multicolumn{1}{|l|}{$O^{*}(4n)$} & \multicolumn{1}{|l|}{$Herm(n,{\Bbb H})$}
& \multicolumn{1}{|l|}{$Sp(1)\times O^{*}(4n-4)$} & \multicolumn{1}{|l|}{$%
O^{*}(4n-4)$} & $4(n-1)$ \\ \hline
\multicolumn{1}{|l|}{$O(2,j)$} & \multicolumn{1}{|l|}{${\Bbb R}^{1,j-1}$} & 
\multicolumn{1}{|l|}{$SL(2,{\Bbb R})\times O(j-2)$} & \multicolumn{1}{|l|}{$%
SL(2,{\Bbb R)}$} & $j-2$ \\ \hline
\multicolumn{1}{|l|}{$E_{7(-25)}$} & \multicolumn{1}{|l|}{$Herm(3,{\Bbb O})$}
& \multicolumn{1}{|l|}{$SO(2,10)$} & \multicolumn{1}{|l|}{$SO(2,10)$} & $16$
\\ \hline
\end{tabular}
}

\subsection{Orbits and representations}

The orbits of the natural action of $L$ on ${\frak n}^{*}=N^{*}$ are
parametrized by pairs of non-negative integers $p^{+},p^{-}$ with $%
p^{+}+p^{-}\leq n$ \cite[2.1]{sahi-expl}. The simplest example of this
parametrization can be observed for $G=Sp(2n,{\Bbb R})$, when the group $N$
can be identified with the vector space of $n\times n$ real symmetric
matrices, and an arbitrary orbit ${\cal O}$ of $L=GL(n,{\Bbb R})$ on $%
N^{*}\simeq N$ is defined by the signature of the symmetric matrix $\xi \in 
{\cal O}$. We write $p$ for the pair $(p^{+},p^{-})$, and ${\cal O}_p$ for
the corresponding orbit. The rank of the orbit is $p^{+}+p^{-},$ which we
denote by $|p|$. If $|p|=n$, the orbit ${\cal O}_p$ is open in $N^{*}$,
otherwise we get small (singular) orbits. By $S_p$ we denote the stabilizer
of $\xi _p\in {\cal O}_p$ in $L$.

Suppose now $|p|<n$. Then $O_p=L/S_p$ has an $L$--equivariant measure $d\mu
_p$ which transforms by some positive character $\delta _p$ of $L$. The main
result of \cite{sahi-expl} associates with each nonzero singular orbit $%
{\cal O}_p$ a unitary irreducible representation $\pi _p$ of $\overline{G}$.
Here $\overline{G}=G$ unless $G=Sp(2n,{\Bbb R)}$ or $O(2,j)$ with $j$ odd
and $\overline{G}$ is a two-fold cover of $G$ in these two cases (for $%
G=Sp(2n,{\Bbb R)}$ we can take $\overline{G}$ to be a real metaplectic
group).

If $H$ is a subgroup of $G$, we write $\overline{H}$ for the inverse image
of $H$ in $\overline{G}$.

The representation $\pi _p$ acts on the Hilbert space $L^2({\cal O}_p,d\mu
_p)$, and actions of the elements of the maximal parabolic subgroup $%
\overline{P}$ can be written in a particularly simple manner -- the action
of the reductive part $\overline{L}$ comes from the action of $\overline{L}$
on ${\cal O}_p$ and the unipotent radical $N$ acts by characters 
\begin{eqnarray}
&&[\pi _p(n)h](\xi _p)=\chi _{\xi _p}(n)\,h(\xi _p),\,n\in N,\xi _p\in {\cal %
O}_p  \label{-pi-on-p} \\
&&[\pi _p(l)h](\xi _p)=\nu _p(l)\delta _p(l)^{-1/2}h(l^{-1}\xi _p),\;l\in 
\overline{L},\xi _p\in {\cal O}_p\;.  \nonumber
\end{eqnarray}
Here $\chi _{\xi _p}$ is the unitary character of the vector space $N$
defined by $\xi _p\in N^{*}$ and $\nu _p$ is a unitary character of $%
\overline{L}$ (trivial on the identity component of $\overline{L}$)\footnote{%
In \cite{sahi-expl} the characters $\nu _p$ and $\delta _p$ are denoted by $%
\mu $ and $\nu $ respectively.}.

\section{Tensor products $\pi _{p_1}\otimes \ldots \otimes \pi _{p_s}$}

We pick $s$ singular orbits ${\cal O}_{p_1},\ldots ,{\cal O}_{p_s}$ such
that $|p_1|+\ldots +|p_s|\leq n$ and consider the tensor product of
associated representations 
\[
\Pi =\bigotimes_{i=1}^s\pi _{p_i}. 
\]
The group $\overline{L}$ acts on the set ${\cal O}^{\prime }\stackrel{\text{%
def}}{=}{\cal O}_{p_1}\times \ldots \times {\cal O}_{p_s}$, and up to a set
of measure zero, ${\cal O}^{\prime }$ is a single $\overline{L}$--orbit.
Note that the set 
\[
{\cal O}\stackrel{\text{def}}{=}{\cal O}_{p_1}+\ldots +{\cal O}%
_{p_s}=\{\zeta \in N^{*}\;|\;\zeta =\tsum\nolimits_{i=1}^s\zeta _{p_i}\text{%
, }\zeta _{p_i}\in {\cal O}_{p_i}\} 
\]
also contains a dense $\overline{L}$--orbit. The representation $\Pi $ acts
in $\bigotimes_{i=1}^sL^2({\cal O}_{p_i},d\mu _{p_i})$, and we can identify
this space with $L^2({\cal O}^{\prime },d\mu ^{\prime })$ where $d\mu
^{\prime }$ is the product measure. If we fix a generic representative $\xi
^{\prime }=(\xi _{p_1},\ldots ,\xi _{p_s})\in {\cal O}^{\prime }$ and set 
\[
\xi =\xi _{p_1}+\ldots +\xi _{p_s}\in {\cal O} 
\]
and $\delta =\prod_{i=1}^s\delta _{p_i}$, $\nu =\prod_{i=1}^s\nu _{p_i}$, we
have the following formulas for the actions of $\Pi |_{\overline{P}}$ on $%
L^2({\cal O}^{\prime },d\mu ^{\prime })$ 
\begin{eqnarray}
&&\Pi (l_0)f(l\xi ^{\prime })=\nu (l_0)\delta (l_0)^{-1/2}h(l_0^{-1}l\xi
^{\prime }),\;l_0\in \overline{L}\text{ }  \label{-pi-prime} \\
&&\Pi (n_0)f(l\xi ^{\prime })=\chi _{l\xi }(n_0)h(l\xi ^{\prime }),\,n_0\in N%
\text{ .}  \nonumber
\end{eqnarray}
Let now $S^{\prime }$ and $S$ be the isotropy subgroups of $\xi ^{\prime }$
and $\xi $, respectively, with respect to the action of $\overline{L}$ on $%
{\cal O}^{\prime }$ and ${\cal O}$. If $|p_1|+\ldots +|p_s|=n,$ the groups $%
S^{\prime }$ and $S$ are reductive. \medskip\ 

{\em Example}. Take $G=U(n,n)$, $s=2$ and $p_1=(k,0)$, $p_2=(0,n-k)$. Then
we can choose $\xi _{p_1}=\left( 
\begin{array}{cc}
I_k & 0 \\ 
0 & 0
\end{array}
\right) $, $\xi _{p_2}=\left( 
\begin{array}{cc}
0 & 0 \\ 
0 & -I_{n-k}
\end{array}
\right) $ and $\xi =\left( 
\begin{array}{cc}
I_k & 0 \\ 
0 & -I_{n-k}
\end{array}
\right) $. It is easy to see that $S=U(k,n-k)$ and $S^{\prime }=U(k)\times
U(n-k)$. The quotient $S/S^{\prime }$ is a Riemannian symmetric space.%
\medskip\ 

\begin{lemma}
$\Pi |_{\overline{P}}\simeq \nu \otimes \limfunc{Ind}_{S^{\prime }N}^{%
\overline{P}}(1\otimes \chi _\xi )$ {\em (}$L^2$-{\em induction)}. Here $\nu 
$ is a character of $\overline{L}$ extended trivially to $\overline{P}$.%
\label{-isom1}
\end{lemma}

\TeXButton{Proof}{\proof}We denote the induced representation $\nu \otimes 
\limfunc{Ind}_{S^{\prime }N}^{\overline{P}}(1\otimes \chi _\xi )$ by $\Pi
^{\prime }$. Then by the definition of the induced representation, $\Pi
^{\prime }$ acts on the space ${\cal C}$ of square-summable functions
satisfying a standard invariance condition 
\begin{equation}
{\cal C}=\{f:\overline{L}N\rightarrow {\Bbb C}\mid \,f(ps^{\prime }n)=\chi
_\xi (n)f(p)\text{ for }p\in \overline{P},s^{\prime }\in S^{\prime },\,n\in
N\}\text{.}  \label{-equivar}
\end{equation}
Let $|p_1|+\ldots +|p_s|<n$ (strict inequality). Then the quasi-invariant
measure on the quotient space $\overline{L}N/S^{\prime }N\simeq \overline{L}%
/S^{\prime }$ is transformed by the character $\delta =\prod_{i=1}^s\delta
_{p_i}$ of $\overline{L}$, and we get 
\begin{eqnarray}
&&\Pi ^{\prime }(l_0)f(ln)=\nu (l_0)\delta (l_0)^{-1/2}f(l_0^{-1}ln)
\label{-pi-2prime} \\
&&\Pi ^{\prime }(n_0)f(ln)=f(n_0^{-1}ln)=f(l(l^{-1}n_0^{-1}ln))\text{ .} 
\nonumber
\end{eqnarray}
We can now define a unitary operator $\Psi :L^2({\cal O}^{\prime },d\mu
^{\prime })\rightarrow {\cal C}$ by setting 
\[
\lbrack \Psi h](ln)=\chi _\xi (n)^{-1}h(l\xi ^{\prime })\text{.}
\]
This operator provides an isometry between $L^2({\cal O}^{\prime },d\mu
^{\prime })$ and ${\cal C}$, and it easy to check that $\Psi $ intertwines
the actions $\Pi $ and $\Pi ^{\prime }$. Indeed, for the actions of $l_0$
this is immediate by inspection of formulas (\ref{-pi-prime}) and (\ref
{-pi-2prime}), and for the actions of $n_0$ we get 
\begin{eqnarray*}
\Pi ^{\prime }(n_0)[\Psi h](ln) &=&[\Psi h](l(l^{-1}n_0^{-1}ln)=\chi _\xi
(l^{-1}n_0^{-1}l)^{-1}\chi _\xi (n)^{-1}h(l\xi ^{\prime }) \\
&=&\chi _\xi (n)^{-1}\chi _{l\xi }(n_0)h(l\xi ^{\prime })=[\Psi \Pi
(n_0)h](ln)\text{. }
\end{eqnarray*}
Computations for $|p_1|+\ldots +|p_s|=n$ are almost identical. In this case
the space $\overline{L}/S^{\prime }$ possesses an $\overline{L}$--invariant
measure, the action of $l_0\in \overline{L}$ is given by 
\[
\Pi ^{\prime }(l_0)f(ln)=f(l_0^{-1}ln)\text{\thinspace ,}
\]
and we set 
\[
\lbrack \Psi h](ln)=\chi _\xi (n)^{-1}\delta (l)^{-1/2}(l\xi ^{\prime })%
\text{.}
\]
A straightforward computation shows that this operator intertwines the
actions of $\Pi $ and $\Pi ^{\prime }$.\TeXButton{End Proof}{\endproof}%
\smallskip\ 

Denote by $\gamma $ a $S^{\prime }$--quasi-regular representation of $S$ in $%
L^2(S/S^{\prime })$. That is 
\begin{equation}
(\gamma (z)f)(x)=f(z^{-1}x)\text{ for }z\in S,\,x\in Y\stackrel{\text{def}}{=%
}S/S^{\prime },f\in L^2(Y).  \label{-sigma}
\end{equation}
Of course, $\gamma =Ind_{S^{\prime }}^S1$, and combining the induction in
stages with the fact that the character $\chi _\xi \,$of $N$ is $SN$-fixed,
we get 
\begin{equation}
\Pi |_{\overline{P}}\simeq \nu \otimes \limfunc{Ind}\nolimits_{S^{\prime
}N}^{\overline{P}}(1\otimes \chi _\xi )=\nu \otimes \limfunc{Ind}%
\nolimits_{SN}^{\overline{P}}((\limfunc{Ind}\nolimits_{S^{\prime
}}^S1)\otimes \chi _\xi )=\nu \otimes \limfunc{Ind}\nolimits_{SN}^{\overline{%
P}}(\gamma \otimes \chi _\xi )\text{.}  \label{-induced-s}
\end{equation}
The groups $S$ and $S^{\prime }$ are, generally speaking, not reductive
(except when $|p_1|+\ldots +|p_s|=n)$. As was discussed in \cite[2.1]
{sahi-expl}, the Lie algebras ${\frak s}$ and ${\frak s}^{\prime }$ of $S$
and $S^{\prime }$, respectively, can be written as 
\begin{eqnarray*}
{\frak s} &=&({\frak l}_1+{\frak g}^{\prime })+{\frak u} \\
{\frak s}^{\prime } &=&({\frak l}_1+{\frak h}^{\prime })+{\frak u\,}\text{,}
\end{eqnarray*}
where ${\frak l}_1,{\frak g}^{\prime },{\frak h}^{\prime }$ are some
reductive Lie algebras, ${\frak h}^{\prime }\subset {\frak g}^{\prime }$ and 
${\frak u}$ is a nilpotent radical common for both ${\frak s}$ and ${\frak s}%
^{\prime }$. Let $G^{\prime }$ and $H^{\prime }$ be the corresponding Lie
groups.

In particular, $X=G^{\prime }/H^{\prime }$ is a reductive homogeneous space,
and we can consider an $H^{\prime }$--quasi-regular representation of $%
G^{\prime }$ on $L^2(X)$ (denoted by $\gamma ^{\prime })$. Then the
representation $\gamma $ of $S$ given by the formula (\ref{-sigma}) can be
obtained by extending $\gamma ^{\prime }$ trivially from $G^{\prime }$ to $S$%
. Now let 
\[
\gamma ^{\prime }\simeq \int_{\widehat{G}^{\prime }}^{\oplus }m(\pi )\pi
\,d\mu (\pi ) 
\]
be a decomposition of a quasi-regular representation $\gamma ^{\prime }$
into a direct integral of unitary irreducible representations of $G^{\prime
} $, where $m:\widehat{G}^{\prime }\rightarrow {\Bbb Z}_{+}$ is a
multiplicity function and $d\mu $ a Plancherel measure for a symmetric space 
$X$. Each irreducible representation $\pi \in \widehat{G}^{\prime }$ can be
extended to an irreducible representation $\pi ^{\vee }$ of $S$. This gives 
\[
\gamma \simeq \int_{\widehat{G}^{\prime }}^{\oplus }m(\pi )\pi ^{\vee
}\,d\mu (\pi ) 
\]
and substituting this into (\ref{-induced-s}), we obtain the decomposition
of Theorem \ref{-th01} 
\begin{equation}
\Pi |_{\overline{P}}\simeq \int_{\widehat{G}^{\prime }}^{\oplus }m(\pi
)\Theta (\pi )\,d\mu (\pi ),  \label{-tensor-decomp}
\end{equation}
where $\Theta (\pi )=\nu \otimes \limfunc{Ind}\nolimits_{SN}^{\overline{P}%
}(\pi ^{\vee }\otimes \chi _\xi )$. Note that representations $\pi $ present
in the formula (\ref{-tensor-decomp}) (i.e., those with $m(\pi )>0$) are $%
H^{\prime }$--spherical representations of $G^{\prime }$.

Mackey theory guarantees that all representations $\Theta (\pi )$ are
unitary irreducible representations of $\overline{P}$ and $\Theta (\pi
)\simeq \Theta (\sigma )$ if and only if $\pi \simeq \sigma $.

The special case $s=2$ deserves some special attention. In this situation $%
\Pi =\pi _p\otimes \pi _q$, where $p=(p^{+},p^{-}),q=(q^{+},q^{-})$. We will
write $G_{pq}^{\prime }$, $H_{pq}^{\prime }$ and $X_{pq}$ for $G^{\prime }$, 
$H^{\prime }$ and $X$, respectively. The quotient space $X_{pq}=G_{pq}^{%
\prime }/H_{pq}^{\prime }$ is then a reductive symmetric space in the sense
of \cite{flensted}. The table below lists these symmetric spaces for
different combinations of $G$, $p$ and $q$ $(\,|p|+|q|\leq n)$ (see 
\cite[16.7]{adams} for the detailed computations in the case of $G=E_7$).%
\smallskip\ 

{\small 
\begin{tabular}{|c|c|c|c|}
\hline
$G$ & $p$ & $q$ & $X_{pq}$ \\ \hline
\multicolumn{1}{|l|}{$Sp(2n,{\Bbb R)}$} & $p$ & $q$ & \multicolumn{1}{|l|}{$%
O(p^{+}+q^{+},p^{-}+q^{-})/[O(p^{+},p^{-})\times O(q^{+},q^{-})]$} \\ \hline
\multicolumn{1}{|l|}{$U(n,n)$} & $p$ & $q$ & \multicolumn{1}{|l|}{$%
U(p^{+}+q^{+},p^{-}+q^{-})/[U(p^{+},p^{-})\times U(q^{+},q^{-})]$} \\ \hline
\multicolumn{1}{|l|}{$O^{*}(4n)$} & $p$ & $q$ & \multicolumn{1}{|l|}{$%
Sp(p^{+}+q^{+},p^{-}+q^{-})/[Sp(p^{+},p^{-})\times Sp(q^{+},q^{-})]$} \\ 
\hline
\multicolumn{1}{|l|}{$O(2,j)$} & $(1,0)$ & $(1,0)$ & \multicolumn{1}{|l|}{$%
SO(j-1)/SO(j-2)$} \\ \hline
\multicolumn{1}{|l|}{} & $(1,0)$ & $(0,1)$ & \multicolumn{1}{|l|}{$%
SO_0(1,j-2)/SO(j-2)$} \\ \hline
\multicolumn{1}{|l|}{$E_{7(-25)}$} & $(1,0)$ & $(1,0)$ & 
\multicolumn{1}{|l|}{$SO(9)/SO(8)$} \\ \hline
\multicolumn{1}{|l|}{} & $(1,0)$ & $(0,1)$ & \multicolumn{1}{|l|}{$%
SO_0(1,8)/SO(8)$} \\ \hline
\multicolumn{1}{|l|}{} & $(1,0)$ & $(2,0)$ & \multicolumn{1}{|l|}{$%
F_{4(-52)}/Spin(9)$} \\ \hline
\multicolumn{1}{|l|}{} & $(1,0)$ & $(0,2)$ & \multicolumn{1}{|l|}{$%
F_{4(-20)}/Spin(9)$} \\ \hline
\multicolumn{1}{|l|}{} & $(1,0)$ & $(1,1)$ & \multicolumn{1}{|l|}{$%
F_{4(-20)}/Spin(1,8)$} \\ \hline
\end{tabular}
} \medskip\ 

\section{Extending $\Theta (\pi )$ to $\overline{G}$}

\subsection{The $N$-spectrum}

In this section we study low-rank representations of $\overline{G}$. For the
classical groups (cases I1, I2, I3 in our list) a complete theory of
low-rank representations can be found in Li's paper \cite{li-lowrank}. We
rely heavily on the ideas and methods of this paper. Our objective here is
to extend the low-rank theory of Li so it can be applied to representations
of the groups $O(2,j)$ and $E_{7(-25)}$.

Consider the restriction of the representation $\Theta (\pi )=\nu \otimes 
\limfunc{Ind}\nolimits_{SN}^{\overline{P}}(\pi ^{\vee }\otimes \chi _\xi )$
to $N$. This restriction decomposes into a direct integral of unitary
characters, and the decomposition is determined by a projection-valued
measure on\thinspace $\widehat{N}=N^{*}$. This measure is supported on the
set ${\cal O}\subset N^{*}$ (an $L$--orbit of $\xi $).

Similarly, for any unitary representation $\tau $ of $\overline{G}$, we can
consider its restriction to the abelian subgroup $N$ and the associated
measure $\mu _\tau $ on $N^{*}$. If $\mu _\tau $ is supported on the single
orbit ${\cal O}_r\subset N^{*}$, we say that $\tau $ is of {\em signature} $%
r=(r^{+},r^{-})$ and write 
\[
\limfunc{sign}\nolimits_N\tau =r\text{.} 
\]
The number $|r|=r^{+}+r^{-}$ is the {\em \ rank }of ${\cal O}_r$. If $\mu
_\tau $ is supported on one or several orbits of rank $k,$ we write $%
\limfunc{rank}\nolimits_N\tau =k$.

It will be convenient to set $\limfunc{sign}t=\left\{ 
\begin{array}{c}
(1,0),t>0 \\ 
(0,1),t<0
\end{array}
\right. $.

{\em Remark}. For the representations of classical groups, the notion of
rank was introduced in \cite{howe-rank} and \cite{li-lowrank}. Our
definition extends it to $G=E_{7(-25)}$. For $G=O(2,j)$ the definition above
differs from the notion of rank in \cite{li-lowrank} due to the different
choice of the parabolic subgroup $P$.

We now take a unitary representation $\sigma $ of $\overline{G}$ and
consider $\sigma |_{\overline{M^{\prime }}N^{\prime }}$. The group $%
N^{\prime }$ is a Heisenberg group defined by an exact sequence 
\[
1\rightarrow ZN^{\prime }\rightarrow N^{\prime }\rightarrow {\Bbb R}%
^{2m}\rightarrow 1\text{,} 
\]
and the multiplication on $N^{\prime }$ defines a standard skew-symmetric
bilinear form on ${\Bbb R}^{2m}$. The group $M^{\prime }$ acts on $N^{\prime
}$ by the automorphisms of $N^{\prime }$, and it also acts trivially on the
center $ZN^{\prime }$. Because of this we can view $M^{\prime }$ as a
subgroup of $Sp(2m,{\Bbb R})$.

Now let $\rho _t$ be a unique representation of $N^{\prime }$ corresponding
in the sense of Stone - von Neumann theorem to the character $\chi _t$ of $%
ZN^{\prime }\simeq {\Bbb R}$, where 
\[
\chi _t(z)=\exp (2\pi itz),\text{ }z\in ZN^{\prime }\text{.} 
\]
We can extend $\rho _t$ to the representation of the semidirect product $%
Sp(2m,{\Bbb R})\widetilde{\;\;\;\cdot }N^{\prime }$ using the corresponding
oscillator representation $\omega _t$ of the metaplectic group $Sp(2m,{\Bbb R%
})\widetilde{\;\;\;}$. This extension restricts to a representation of a
semidirect product $\overline{M^{\prime }}N^{\prime }$, and we denote this
restriction by $\widetilde{\rho }_t$.

By the results of \cite{howe-moore} the subspace of $ZN^{\prime }$--fixed
vectors is invariant under the action of $\sigma (\overline{G}),$ and
without loss of generality we may assume that $\sigma $ has no $ZN^{\prime }$%
--fixed vectors. Then, according to the Mackey theory, $\sigma |_{\overline{%
M^{\prime }}N^{\prime }}$ decomposes into representations of the form $%
\kappa _t\otimes \widetilde{\rho }_t$, where all $\kappa _t$, $t\in {\Bbb R}%
^{*}$ are unitary representations of $\overline{M^{\prime }}$. We can write 
\[
\sigma |_{\overline{G}_{-}\,N^{\prime }}=\int_{{\Bbb R}^{*}}^{\oplus }\kappa
_t\otimes \widetilde{\rho }_t\,dt. 
\]
We will now describe the $N$-spectrum of $\widetilde{\rho }_t$. It is known
that a real vector space $N$ is endowed with a structure of a simple
formally real Jordan algebra with the unit (denoted by $e$), and $L$ is the
structure group of this Jordan algebra. Then $ZN^{\prime }$ is a
one-dimensional subalgebra of $N$ generated by the primitive idempotent $c$.

This idempotent determines the Peirce decomposition of $N$ \cite[IV.I]{fk}: 
\[
N=N(c,1)+N(c,1/2)+N(c,0)\text{.} 
\]
Observe that $N(c,1)=ZN^{\prime }$, $N(c,1/2)={\Bbb R}^m$ and $N(c,0)=N_{-}.$

{\em Example}. Take $G=E_{7(-25)}$. Then $N=Herm(3,{\Bbb O)}$, $c=\left[ 
\begin{array}{ccc}
1 & 0 & 0 \\ 
0 & 0 & 0 \\ 
0 & 0 & 0
\end{array}
\right] $ and the corresponding Peirce decomposition is 
\[
N={\Bbb R}c+{\Bbb O}^2+Herm(2,{\Bbb O)}\text{.} 
\]
Hence $N(c,1/2)={\Bbb O}^2={\Bbb R}^{16}$ and $N(c,0)=Herm(2,{\Bbb O)=R}%
^{1,9}$ and this Jordan algebra is in fact the nilradical $N_{-}$ of the
parabolic subgroup $P_{-}=[{\Bbb R}^{*}\times SO(1,9)]\cdot {\Bbb R}^{1,9}$
inside $G_{-}=SO(2,10)$.\medskip\ 

The action of $\widetilde{\rho _t}$ on $N(c,1)$ and $N(c,1/2)$ is easy to
describe: these spaces lie inside $N^{\prime }$ and we can model an
irreducible representation $\rho _t$ of $N^{\prime }$ on the Hilbert space $%
L^2(N(c,1/2))$. To distinguish between an element $x$ of $N(c,1/2)$ and the
corresponding vector from ${\Bbb R}^m$, we shall write $\widehat{x}$ for the
latter. Then 
\begin{eqnarray}
&&\widetilde{\rho _t}(n_1c)f(x)=\chi _t(n_1)f(x)=\chi _t(tr_N[n_1c])f(x)%
\text{, }n_1\in {\Bbb R}  \label{-rho-t1} \\
&&\widetilde{\rho _t}(n_{1/2})f(x)=\chi _t(\widehat{n}_{1/2}\cdot \widehat{x}%
)f(x)=\chi _t(\tfrac 12tr_N[n_{1/2}x])f(x)\text{, }n_{1/2}\in N(c,1/2). 
\nonumber
\end{eqnarray}
Here $tr_N$ is the standard trace functional on the Jordan algebra $N$.

Now take $z_0\in Sym(m,{\Bbb R})\subset Sp(2m,{\Bbb R})$. The action of the
oscillator representation $\omega _t(z_0)$ on $L^2({\Bbb R}^m)$ is given by
the formula 
\[
\omega _t(z_0)f(x)=\chi _t(\tfrac 12\widehat{x}z_0\widehat{x}^t)f(x). 
\]
Observe that $\widehat{x}z_0\widehat{x}^t=tr_{Sym(m,{\Bbb R})}[(e-c)x^2z_0].$

Recall that $G_{-}\subset Sp(2m,{\Bbb R})$ and 
\[
N(c,0)=N_{-}=P_{-}\cap Sym(m,{\Bbb R})\text{.} 
\]
For $n_0\in N(c,0)$ and $x\in N(c,1/2)$ we have 
\[
tr_{Sym(m,{\Bbb R})}[(e-c)x^2n_0]=tr_{N_{-}}[(e-c)x^2n_0] 
\]
and 
\begin{equation}
\widetilde{\rho _t}(n_0)f(x)=\chi _t(\tfrac 12tr_N[(e-c)x^2n_0])f(x)\text{.}
\label{-rho-t2}
\end{equation}

Combining formulas (\ref{-rho-t1}) and (\ref{-rho-t2}), we can write the
formula for $\widetilde{\rho _t}(n^0)$, where $n^0=n_1c+n_{1/2}+2n_0$ . 
\begin{equation}
\widetilde{\rho _t}(n^0)f(x)=\chi _t(tr_N([c+\tfrac 12x+\tfrac
14(e-c)x^2]n^0))f(x)\text{.}  \label{-rho-tn}
\end{equation}
We can identify $N^{*}$ and $N$ by setting $\phi (n^{\prime })(n^{\prime
\prime })=tr_N(n^{\prime }n^{\prime \prime })$ for $n,n^{\prime }\in N$. It
follows from the formula (\ref{-rho-tn}) that the $N$-spectrum of $%
\widetilde{\rho _t}$ is supported on the elements of the form $n_t(x),\,x\in
N(c,1/2)$, where 
\[
n_t(x)=t(c+\tfrac 12x+\tfrac 14(e-c)x^2)\text{.} 
\]
For an arbitrary element $x^{\prime }\in N(c,1/2)$ there exists a special
element of the structure group $L$,$\,$called the Frobenius transformation
and denoted by $\tau (x^{\prime })$. According to Lemma VI.3.1 of \cite{fk}, 
$\tau (x^{\prime })n_t(x)=n_1^{\prime }+n_{1/2}^{\prime }+n_0^{\prime }$,
where 
\begin{eqnarray*}
&&n_1^{\prime }=tc \\
&&n_{1/2}^{\prime }=t(2x^{\prime }c+\tfrac 12x) \\
&&n_0^{\prime }=t(2(e-c)x^{\prime 2}c+(e-c)x^{\prime }x+\tfrac 14(e-c)x^2)\,.
\end{eqnarray*}
In particular, $\tau (-\frac x2)n_t(x)=tc.$

We can now describe the $N$-spectrum of $\kappa _t\otimes \widetilde{\rho _t}
$. If the $N_{-}$-spectrum of $\kappa _t$ is supported on a set ${\cal O}%
(\kappa _t)$, then the support of the $N$-spectrum of $\kappa _t\otimes 
\widetilde{\rho _t}$ consists of the elements $n_t(x)+n_{-}$ where $x\in
N(c,1/2),n_{-}\in {\cal O}(\kappa _t)$. Then 
\[
\tau (-\tfrac x2)(n_t(x)+n_{-})=tc+n_{-}\;\text{.} 
\]
Suppose now that $\limfunc{sign}\nolimits_N\sigma =r$, where $%
r=(r^{+},r^{-}) $. Then $\limfunc{sign}\nolimits_N(\kappa _t\otimes 
\widetilde{\rho _t})=r$, i.e., 
\begin{equation}
\limfunc{sign}\nolimits_N(tc+n_{-})=r\text{.}  \label{-sign-n}
\end{equation}
It is easy to see that (\ref{-sign-n}) implies $\limfunc{sign}t+\limfunc{sign%
}\nolimits_{N_{-}}n_{-}=r$, i.e., 
\[
\limfunc{sign}\nolimits_{N_{-}}\kappa _t=\left\{ 
\begin{array}{c}
(r^{+}-1,r^{-}),\;t>0 \\ 
(r^{+},r^{-}-1),\;t<0
\end{array}
\right. \text{.} 
\]
We summarize this discussion in the following

\begin{lemma}
Let $\sigma $ be a representation of $\overline{G}$, $\limfunc{sign}%
\nolimits_N\sigma =r$ and $\sigma |_{\overline{G}_{-}\,N^{\prime
}}=\dint\limits_{{\Bbb R}^{*}}^{\oplus }\kappa _t\otimes \widetilde{\rho }%
_t\,dt.$ Then for any $t\in {\Bbb R}^{*}$ the $N_{-}$-spectrum of the
representation $\kappa _t$ is supported on a single $L_{-}$-orbit, and $%
\limfunc{sign}\nolimits_{N_{-}}\kappa _t=r-\limfunc{sign}t$.\label
{-lemma-integral}
\end{lemma}

\subsection{Von\ Neumann algebras}

Let $\tau $ be a representation of some subgroup $H$ of $\overline{G}$. By $%
{\cal A}(\tau ,H_0)$ we denote the von Neumann algebra generated by the
operators $\tau (h),h\in H_0$, where $H_0$ is a subgroup of $H.$

To proceed further we need

\begin{lemma}
Assume that ${\cal A}(\kappa _t,\overline{G}_{-})={\cal A}(\kappa _t,%
\overline{P}_{-})$ for all $t\in {\Bbb R}^{*}$. Then 
\[
{\cal A}(\dint\limits_{{\Bbb R}^{*}}^{\oplus }\kappa _t\otimes \widetilde{%
\rho }_t\,dt\,,\overline{G}_{-})\subseteq {\cal A}(\dint\limits_{{\Bbb R}%
^{*}}^{\oplus }\kappa _t\otimes \widetilde{\rho }_t\,dt\,,\overline{P}%
_{-}N^{\prime })\text{.}
\]
\label{-lemma-vn}
\end{lemma}

\TeXButton{Proof}{\proof} The representation $\rho _t$ is an irreducible
representation of $N^{\prime }$, therefore ${\cal A}(\rho _t,N_{-})$ is the
full algebra of bounded operators on $L^2({\Bbb R}^m).$ Consider the algebra 
${\cal A}(\kappa _t\otimes \widetilde{\rho }_t,\overline{G}_{-})$. This
algebra is generated by operators 
\begin{equation}
\kappa _t(g_{-})\otimes \widetilde{\rho }_t(g_{-}),g_{-}\in \overline{G}_{-}.
\label{-nuxrho}
\end{equation}
All these operators lie inside ${\cal A}(\kappa _t\otimes \widetilde{\rho }%
_t,\overline{P}_{-}N^{\prime })$. Indeed, the algebra ${\cal A}(\kappa
_t\otimes \widetilde{\rho }_t,\overline{P}_{-}N^{\prime })$ contains the set 
${\cal B}$ of all operators of the form $\kappa _t(p^{-})\otimes a$, where $a
$ is an arbitrary bounded operator on $L^2({\Bbb R}^m)$ and $p^{-}\in 
\overline{P}_{-}$. Combining this fact with the assumption ${\cal A}(\kappa
_t,\overline{G}_{-})={\cal A}(\kappa _t,\overline{P}_{-})$, we conclude that
the von Neumann{\em \ }algebra generated by ${\cal B}$ already contains all
operators (\ref{-nuxrho}).

Hence ${\cal A}(\kappa _t\otimes \widetilde{\rho }_t,\overline{G}%
_{-})\subseteq {\cal A}(\kappa _t\otimes \widetilde{\rho }_t,\overline{P}%
_{-}N^{\prime })$ and 
\[
{\cal A}(\dint\limits_{{\Bbb R}^{*}}^{\oplus }\kappa _t\otimes \widetilde{%
\rho }_t\,dt\,,\overline{G}_{-})\subseteq \dint\limits_{{\Bbb R}%
^{*}}^{\oplus }{\cal A}(\kappa _t\otimes \widetilde{\rho }_t\,,\overline{G}%
_{-})\,dt\subseteq \dint\limits_{{\Bbb R}^{*}}^{\oplus }{\cal A}(\kappa
_t\otimes \widetilde{\rho }_t\,,\overline{P}_{-}N^{\prime })\,dt\;. 
\]
But the representations $\widetilde{\rho }_t$ are irreducible and
nonisomorphic for different $t$, and 
\[
\dint\limits_{{\Bbb R}^{*}}^{\oplus }{\cal A}(\kappa _t\otimes \widetilde{%
\rho }_t\,,\overline{P}_{-}N^{\prime })\,dt={\cal A}(\dint\limits_{{\Bbb R}%
^{*}}^{\oplus }\kappa _t\otimes \widetilde{\rho }_t\,dt\,,\overline{P}%
_{-}N^{\prime })\text{. \ \TeXButton{End Proof}{\endproof}} 
\]

Observe that $\overline{P}_{-}N^{\prime }$ is a subgroup of $\overline{P}$.

The next theorem is an analogue of \cite[4.3]{li-lowrank}.

\begin{theorem}
Let $\sigma $ be a representation of $\overline{G}$ , $\limfunc{sign}%
\nolimits_N\sigma =r$ and $0<\left| r\right| <n$ (i.e., $\sigma $ is a
low-rank representation of $\overline{G}$)$.$ Then ${\cal A}(\sigma ,%
\overline{G})={\cal A}(\sigma ,\overline{P})$ $.$\label{-theorem-vn}
\end{theorem}

\TeXButton{Proof}{\proof}The groups $\overline{G}_{-}$ and $\overline{P}$
together generate $\overline{G}$, and it suffices to check that 
\begin{equation}
{\cal A}(\sigma ,\overline{G}_{-})\subseteq {\cal A}(\sigma ,\overline{P})%
\text{.}  \label{-vn-inclusion}
\end{equation}
But $\sigma |_{\overline{G}_{-}\,N^{\prime }}=\dint\limits_{{\Bbb R}%
^{*}}^{\oplus }\kappa _t\otimes \widetilde{\rho }_t\,dt$, and by Lemma \ref
{-lemma-vn}, the assertion (\ref{-vn-inclusion}) follows immediately if we
can show that ${\cal A}(\kappa _t,\overline{G}_{-})={\cal A}(\kappa _t,%
\overline{P}_{-})$ for all $t\in {\Bbb R}^{*}$. By Lemma \ref
{-lemma-integral} all $\kappa _t$ are representations of rank $\left|
r\right| -1$ of the group $\overline{G}_{-}$, and we can apply the same line
of reasoning to them.

After $\left| r\right| $ steps of this process, we reduce the statement of
the theorem to the case of representations of rank 0 for some group $%
\overline{G}_0$, where $G_0$ is belongs to one of the families I1-I4. Any
representation $\tau $ of rank 0 decomposes over characters of $\overline{G}%
_0$ \cite{howe-moore} and it is well known that any character of $\overline{%
G_0}$ is determined by its restriction to the Siegel parabolic $\overline{P}%
_0$ (e.g., \cite[4.2]{li-lowrank}). Therefore, ${\cal A}(\tau ,\overline{G}%
_0)={\cal A}(\tau ,\overline{P}_0)$.\ \ \TeXButton{End Proof}{\endproof}%
\medskip\ 

We now return to the problem of decomposing representation $\Pi
=\bigotimes_{i=1}^s\pi _{p_i}$. The restriction of this representation on $%
\overline{P}$ is given by (\ref{-tensor-decomp}), and for any $\Theta (\pi
)=\nu \otimes \limfunc{Ind}\nolimits_{SN}^{\overline{P}}(\pi ^{\vee }\otimes
\chi _\xi )$ in the decomposition (\ref{-tensor-decomp}) 
\[
\limfunc{sign}\nolimits_N\Theta (\pi )=\limfunc{sign}\nolimits_N\xi
=\tsum\nolimits_{i=1}^sp_i\,. 
\]
Therefore $\Pi $ can be decomposed over the irreducible representations of $%
\overline{G}$ of signature $\tsum\nolimits_{i=1}^sp_i$.

Assume $\sum_{i=1}^s\left| p_i\right| <n$. Then by Theorem \ref{-theorem-vn}%
, any two non-isomorphic irreducible representation from the spectrum of $%
\Pi $ restrict to non-isomorphic irreducible representations of $\overline{P}
$. Therefore the $\overline{P}$-decomposition (\ref{-tensor-decomp}) gives
rise to a $\overline{G}$-decomposition 
\begin{equation}
\Pi \simeq \int_{\widehat{G}^{\prime }}^{\oplus }m(\pi )\theta (\pi )\,d\mu
(\pi ),  \label{-tensor-g}
\end{equation}
where $\theta (\pi )$ is defined for almost every $\pi $ (with respect to $%
d\mu )\,$as a unique irreducible representation of $\overline{G}$ determined
by the condition $\theta (\pi )|_{\overline{P}}=\nu \otimes \limfunc{Ind}%
\nolimits_{SN}^{\overline{P}}(\pi ^{\vee }\otimes \chi _\xi )$. Obviously $%
\theta (\pi )\simeq \theta (\sigma )$ if and only if $\pi \simeq \sigma $.

\section{Representations of maximal rank}

The statement of Theorem \ref{-theorem-vn} is certainly false for the
representations of maximal possible rank, i.e., when $\limfunc{sign}%
\nolimits_N\sigma =r$ and $\left| r\right| =n$. Nevertheless, a $\overline{G}
$-decomposition (\ref{-tensor-g}) can be constructed even when $%
\sum_{i=1}^s\left| p_i\right| =n$.

Consider $\sigma =\sigma _1\otimes \sigma _2$, where $\sigma _1$ and $\sigma
_2$ are representations of $\overline{G}$, $\limfunc{sign}\nolimits_N\sigma
_1=r_1=(r^{+},r^{-})$, $\left| r_1\right| =n-1$ and $\limfunc{sign}%
\nolimits_N\sigma _2=(1,0)$. Then $\sigma _1|_{\overline{M}^{\prime
}\,N^{\prime }}=\dint\limits_{{\Bbb R}^{*}}^{\oplus }\kappa _t\otimes 
\widetilde{\rho }_t\,dt$ and $\sigma _2|_{\overline{M}^{\prime }\,N^{\prime
}}=\dint\limits_{{\Bbb R}_{+}^{*}}^{\oplus }\kappa _u^{\prime }\otimes 
\widetilde{\rho }_u\,du$, and 
\begin{equation}
\sigma |_{\overline{M}^{\prime }\,N^{\prime }}=\iint_{{\Bbb R}^{*}\times 
{\Bbb R}_{+}^{*}}^{\oplus }(\kappa _t\otimes \kappa _u^{\prime })\otimes (%
\widetilde{\rho }_t\otimes \widetilde{\rho }_u)\,dt\,du\text{.}
\label{-doubleint}
\end{equation}
For $t+u\neq 0$ we have $(\widetilde{\rho }_t\otimes \widetilde{\rho }%
_u)|_{N^{\prime }}=\rho _t\otimes \rho _u\simeq 1\otimes \rho _{t+u}$, where 
$1$ is a trivial representation of $N^{\prime }$ on $L^2({\Bbb R}^m)$.

$(\widetilde{\rho }_t\otimes \widetilde{\rho }_u)|_{\overline{M}^{\prime
}}=\omega _t\otimes \omega _u$ $\simeq \omega _{t,u}^{\prime }\otimes \omega
_{t,u}^{\prime \prime }$ where $\omega _{t,u}^{\prime \prime }=\left\{ 
\begin{array}{c}
\omega _{+},t+u>0 \\ 
\omega _{-},t+u<0
\end{array}
\right. $ and $\omega _{t,u}^{\prime }=\left\{ 
\begin{array}{c}
\omega _{+},tu/(t+u)>0 \\ 
\omega _{-},tu/(t+u)<0
\end{array}
\right. $. Here $\omega _{+}$ and $\omega _{-}$ are the restrictions of two
nonisomorphic oscillator representations of $Sp(2m,{\Bbb R})\widetilde{\;\;\;%
}$ to $\overline{M}^{\prime }$.

Then $\widetilde{\rho }_t\otimes \widetilde{\rho }_u\simeq \tau
_{t,u}\otimes \widetilde{\rho }_{t+u}$, where $\tau _{t,u}(N^{\prime })$
acts trivially on $L^2({\Bbb R}^m)$ and $\tau _{t,u}(\overline{M}^{\prime })$
acts by $\omega _{t,u}^{\prime }$. The set $t+u=0$ has measure 0 in ${\Bbb R}%
^{*}\times {\Bbb R}^{*}$ and after a change of variables $t+u=v$ the
decomposition (\ref{-doubleint}) becomes 
\[
\sigma |_{\overline{M}^{\prime }N^{\prime }}=\iint_{{\cal D}}^{\oplus
}(\kappa _t\otimes \kappa _{v-t}^{\prime }\otimes \tau _{t,v-t})\otimes 
\widetilde{\rho }_v\,dt\,dv\text{,} 
\]
where ${\cal D}=\{(t,v)\;|\,t\neq 0,v\neq 0,\,v>t\}$.

If we set $\lambda _v=\dint\limits_{(-\infty ,v)}^{\oplus }\kappa _t\otimes
\kappa _{v-t}^{\prime }\otimes \tau _{t,v-t}\,dt$, the preceding formula can
be rewritten as 
\begin{equation}
\sigma |_{\overline{M}^{\prime }N^{\prime }}=\dint\limits_{{\Bbb R}%
^{*}}^{\oplus }\lambda _v\otimes \widetilde{\rho }_v\,dv.  \label{-dirint-l}
\end{equation}

By Lemma \ref{-lemma-integral} all representations $\kappa _t$ have
signature $r_1-\limfunc{sign}t$, and all $\kappa _{v-t}^{\prime }$ are of
rank 0, i.e., decomposable over characters. Therefore 
\begin{equation}
\lambda _v|_{\overline{M}^{\prime }}=\left\{ 
\begin{array}{l}
\kappa _v^{-}\otimes \omega _{+},\limfunc{sign}_{N_{-}}\kappa
_v^{-}=(r^{+},r^{-}-1)\text{ if }v<0 \\ 
(\kappa _v^{-}\otimes \omega _{-})\oplus (\kappa _v^{+}\otimes \omega _{+}),
\\ 
\;\;\;\;\;\;\;\;\limfunc{sign}_{N_{-}}\kappa _v^{-}=(r^{+},r^{-}-1),\limfunc{%
sign}_{N_{-}}\kappa _v^{+}=(r^{+}-1,r^{-})\text{ if }v>0.
\end{array}
\right.  \label{-lambda}
\end{equation}

{\em Remark}. If the signature $r$ is semi-definite (i.e. $r=(\left|
r\right| ,0)$ or $(0,\left| r\right| )$), some of the signatures in the
formula above will involve negative numbers, which is of course impossible.
To simplify notation, we agree that in this case corresponding summands are
simply absent from the decomposition (\ref{-dirint-l}).

\begin{lemma}
Let $\sigma =\sigma _1\otimes \sigma _2$, where $\sigma _1$ and $\sigma _2$
are representations of $\overline{G}=\overline{O(2,j)}$, $\limfunc{sign}%
_N\sigma _1=r$, $|r|=1$ and $\limfunc{sign}_N\sigma _2=(1,0)$. Then ${\cal A}%
(\sigma ,\overline{G})={\cal A}(\sigma ,\overline{P})$.\label{-l-so}
\end{lemma}

\TeXButton{Proof}{\proof}Consider $P^1=P_{-}\times O(j-2)$ -- a parabolic
subgroup of $M^{\prime }=SL(2)\times O(j-2)$. It suffices to prove that for
all $v$ 
\[
{\cal A}(\lambda _v,\overline{M}^{\prime })={\cal A}(\lambda _v,\overline{P^1%
})\text{.}
\]
Indeed, this fact combined with formula (\ref{-dirint-l}) and Lemma \ref
{-lemma-vn} gives ${\cal A}(\sigma ,\overline{M}^{\prime })\subseteq {\cal A}%
(\sigma ,\overline{P})$, and the statement of the lemma follows.

Analysis of (\ref{-lambda}) shows that $\lambda _v|_{\overline{M}^{\prime
}}=\chi _v\otimes \omega ,$ where $\chi _v$ decomposes over characters and $%
\omega $ is an oscillator representation restricted to $\overline{M^{\prime }%
}$. Without loss of generality we may take $\omega =\omega _{+}$. Two
factors of $M^{\prime }$ form a dual reductive pair inside $Sp(2(j-2),{\Bbb R%
})$ and the spectrum of $\omega _{+}$ is very well known: $\omega
_{+}=\bigoplus_i$ $\eta _1^{(i)}\otimes \eta _2^{(i)}$, where $\eta _1^{(i)}$
and $\eta _2^{(i)}$ are irreducible highest weight representations of $%
\overline{SL(2)}$ and $O(j-2)$ respectively, and each $\eta _1^{(i)}$ and $%
\eta _2^{(i)}$ occurs only once in the decomposition. Observe that each $%
\eta _1^{(i)}|_{\overline{P}_{-}}$ is irreducible. Therefore ${\cal A}(\eta
_1^{(i)}\otimes \eta _2^{(i)},\overline{M}^{\prime })={\cal A}(\eta
_1^{(i)}\otimes \eta _2^{(i)},\overline{P^1})$, and ${\cal A}(\omega ,%
\overline{M}^{\prime })={\cal A}(\omega ,\overline{P^1})$. Similarly, ${\cal %
A}(\chi _v\otimes \omega ,\overline{M}^{\prime })={\cal A}(\chi _v\otimes
\omega ,\overline{P^1})$. Hence any irreducible component of $\chi _v\otimes
\omega $ is irreducible when restricted to $\overline{P^1}$ and uniquely
determined by this restriction, and ${\cal A}(\chi _v\otimes \omega ,%
\overline{M}^{\prime })={\cal A}(\chi _v\otimes \omega ,\overline{P^1})$. 
\TeXButton{End Proof}{\endproof}\medskip\ 

{\em Remark}. It is easy to see (by inspection of the above argument) that
the statement of the lemma remains true if we replace $\sigma =\sigma
_1\otimes \sigma _2$ with $\sigma =\bigoplus_{i=1}^k\sigma _1^{(i)}\otimes
\sigma _2^{(i)}$, where $\limfunc{sign}_N\sigma _1^{(i)}=r$ and $\limfunc{%
sign}_N\sigma _2^{(i)}=(1,0)$, $1\leq i\leq k$. We can also replace $%
\overline{G}=\overline{O(2,j)}$ with $\overline{SO(2,j)}$.

\begin{lemma}
Let $\sigma =\sigma _1\otimes \sigma _2$, where $\sigma _1$ and $\sigma _2$
are representations of $G=E_{7(-25)}$, $\limfunc{sign}_N\sigma _1=r$, $|r|=2$
and $\limfunc{sign}_N\sigma _2=(1,0)$. Then ${\cal A}(\sigma ,G)={\cal A}%
(\sigma ,P)$.\label{-l-e7}
\end{lemma}

\TeXButton{Proof}{\proof}In this case $M^{\prime }=G_{-}=SO(2,10).$ Once
again, it suffices to check that for all $v\in {\Bbb R}^{*}$%
\begin{equation}
{\cal A}(\lambda _v,G_{-})={\cal A}(\lambda _v,P_{-})\text{.}
\label{-vn-lambda}
\end{equation}

From (\ref{-lambda}) we see that $\lambda _v$ is either a tensor product of
two representations of rank 1 (in this case the assertion of (\ref
{-vn-lambda}) follows immediately from Lemma \ref{-l-so}) or 
\[
\lambda _v=(\kappa _v^{-}\otimes \omega _{-})\oplus (\kappa _v^{+}\otimes
\omega _{+})\text{,} 
\]
where $\limfunc{sign}_{N_{-}}\kappa _v^{-}=(1,0)$, $\limfunc{sign}%
_{N_{-}}\kappa _v^{+}=(0,1)$. This can occur only when $r=(1,1),v>0$. But $%
\limfunc{sign}_{N_{-}}\omega _{-}=(0,1)$, $\limfunc{sign}_{N_{-}}\omega
_{+}=(1,0)$ and we find ourselves in a situation described in the remark to
Lemma \ref{-l-so}.

Therefore, (\ref{-vn-lambda}) holds for all $v$. \TeXButton{End Proof}
{\endproof}\medskip

\begin{corollary}
Let $\sigma =\pi _{p_1}\otimes \pi _{p_2}\otimes \pi _{p_3}$ be a
representation of $G=E_{7(-25)}$, $|p_1|=1$, $|p_2|=1$, $p_3=(1,0)$. Then $%
{\cal A}(\sigma ,G)={\cal A}(\sigma ,P)$.\label{-l-e7-three}
\end{corollary}

\TeXButton{Proof}{\proof}Set $\sigma _1=\pi _{p_1}\otimes \pi _{p_2}$, $%
\sigma _2=\pi _{p_3}$. Then $\limfunc{sign}_N\sigma _1=p_1+p_2$, $\limfunc{%
sign}_N\sigma _2=(1,0)$ and the lemma above can be applied. 
\TeXButton{End Proof}{\endproof}

We now return to our study of the tensor product $\Pi =\bigotimes_{i=1}^s\pi
_{p_i}$, $\sum_{i=1}^s\left| p_i\right| \leq n$.

\begin{theorem}
${\cal A}(\Pi ,\overline{G})={\cal A}(\Pi ,\overline{P})$.
\end{theorem}

\TeXButton{Proof}{\proof}If $\sum_{i=1}^s\left| p_i\right| <n$, the
statement of this theorem follows from Theorem \ref{-theorem-vn}. Hence we
can restrict our attention to the case $\sum_{i=1}^s\left| p_i\right| =n$.

If $\overline{G}=\overline{O(2,j)}$, the only possible case is $s=2$, and we
may always assume $p_2=(1,0)$ and apply Lemma \ref{-l-so}. Similarly, for $%
\overline{G}=E_{7(-25)}$ we can take $p_s=(1,0)$, and the theorem follows
from Lemma \ref{-l-e7} for $s=2$ and Corollary \ref{-l-e7-three} for $s=3$.

Finally, in the classical cases (I1-I3) the statement follows immediately
from \cite[4.7-4.8]{li-lowrank}. Indeed, for these groups each of the
representations $\pi _{p_i}$ appears in the Howe duality correspondence for
an appropriate stable range dual pair $(G_i^{\prime },G)$ and all
irreducible representations from the spectrum of $\Pi $ appear in the
duality correspondence for the pair $(G^{\prime },G)$, which is still in the
stable range. Therefore any irreducible representation from the spectrum of $%
\Pi $ is irreducible when restricted to $\overline{P}$ and uniquely
determined by this restriction, and ${\cal A}(\Pi ,\overline{G})={\cal A}%
(\Pi ,\overline{P})$.\TeXButton{End Proof}{\endproof}\ 

Therefore the $\overline{P}$-decomposition (\ref{-tensor-decomp}) gives rise
to a $\overline{G}$-decomposition of $\Pi $ with respect to the same measure 
$d\mu $ and multiplicity function $m(\pi )$ 
\begin{equation}
\Pi \simeq \int_{\widehat{G}^{\prime }}^{\oplus }m(\pi )\theta (\pi )\,d\mu
(\pi ).  \label{-tensor-g1}
\end{equation}
Comparing (\ref{-tensor-decomp}) and (\ref{-tensor-g1}) we see that $\theta
(\pi )$ is a unitary irreducible representation of $\overline{G}$ which can
be defined (for almost every $\pi $ with respect to $d\mu )\,$as a unique
irreducible representation from the spectrum of $\Pi $ satisfying the
condition $\theta (\pi )|_{\overline{P}}=\Theta (\pi )$, where $\Theta (\pi
)=\nu \otimes \limfunc{Ind}\nolimits_{SN}^{\overline{P}}(\pi ^{\vee }\otimes
\chi _\xi )$.

Theorem \ref{-g-decomp} is thus proved.\smallskip\

\end{document}